\documentclass[final,3p,times]{elsarticle}
\usepackage{amsfonts, amsthm}
\usepackage{latexsym}
\usepackage{amssymb}
\usepackage[mathscr]{eucal}
\usepackage{epsfig}
\usepackage{lscape}
\usepackage{amssymb}
\usepackage{amsmath}
\usepackage{graphicx}
\usepackage{pstricks,pst-plot,pst-node,pst-slpe}
\usepackage{amsfonts}%
\providecommand{\U}[1]{\protect\rule{.1in}{.1in}}

\newtheorem{theorem}{Theorem}
\newtheorem{proposition}[theorem]{Proposition}
\newtheorem{lemma}[theorem]{Lemma}

\def\a{\alpha}

\def\be{\begin{equation}}
\def\ee{\end{equation}}
\def\bea{\begin{eqnarray}}
\def\eea{\end{eqnarray}}
\def\bean{\begin{eqnarray*}}
\def\eean{\end{eqnarray*}}

\def\proof#1. {\par
\ifdim\lastskip<15pt
\removelastskip\penalty-200
\vskip5pt plus3pt minus3pt
\fi
{\def\a{#1}
\ifx\a\empty
{\noindent\bf Proof.}
\else
{\noindent\bf Proof of #1.}
\fi}\enspace}

\journal {Journal of Computational and Applied Mathematics}
\begin{document}
\begin{frontmatter}
\author{D Dominici\fnref{label1}}\ead{dominicd@newpaltz.edu}\fntext[label1]{Research by this author was supported by a Humboldt Research Fellowship for Experienced Researchers from the Alexander von Humboldt Foundation.} \address{Department of Mathematics, State University of New York at New Paltz, 1 Hawk Dr. Suite 9, New Paltz, NY 12561-2443, USA}
\author{SJ Johnston}\ead{johnssj@unisa.ac.za}\address{Department of Mathematical Sciences, University of South Africa, PO Box 392, UNISA 0003, South Africa}
\author{K Jordaan\fnref{label2}}\ead{kerstin.jordaan@up.ac.za}\fntext[label2]{Research by this author was partially supported by the National Research Foundation under grant number 2054423.} \address{Department of Mathematics and Applied Mathematics, University
of Pretoria, Lynnwood Road, Pretoria, 0002, South Africa}
\title{Real zeros of $_2F_1$ hypergeometric polynomials}
\begin{abstract}
We use a method based on the division algorithm to determine all the values of the real parameters $b$ and $c$ for which the hypergeometric polynomials $_2F_1(-n, b; c; z)$ have $n$ real, simple zeros. Furthermore, we use the quasi-orthogonality of Jacobi polynomials to determine the intervals on the real line where the zeros are located.
\end{abstract}
\end{frontmatter}

\section{Introduction}

The ${}_{2}F_{1}$ hypergeometric function is defined by (cf. \cite{aar})
\[
{}_{2}F_{1}\left(  a,b;c;z\right)  =1+\sum_{k=1}^{\infty}\frac{(a)_{k}(b)_{k}%
}{(c)_{k}}\,\frac{z^{k}}{k!},\quad|z|<1,
\]
where $a,\,b$ and $c$ are complex parameters, $-c\notin\mathbb{N}_{0}=\left\{
0,1,2,\ldots\right\}  $ and
\[
(\alpha)_{k}=%
\begin{cases}
\alpha(\alpha+1)\ldots(\alpha+k-1) & \quad,\quad k\in\mathbb{N},\medspace\\
1 & \quad,\quad k=0,\;\alpha\neq0
\end{cases}
\]
is Pochhammer's symbol. This series converges when $|z|<1$ and also when $z=1$
provided that $\text{Re}(c-a-b)>0$ and when $z=-1$ provided that
$\text{Re}(c-a-b+1)>0$. When one of the numerator parameters is equal to a
nonpositive integer, say $a=-n,\;n\in\mathbb{N}_{0}$, the series terminates
and the function is a polynomial of degree $n$ in $z$.

The problem of describing the zeros of the polynomials ${}_{2}F_{1}\left(
-n,b;c;z\right)  $ when $b$ and $c$ are complex arbitrary parameters, has not been
solved. Even when $b$ and $c$ are both real, the only cases that have been
fully analyzed impose additional restrictions on $b$ and $c$. Recent
publications (cf. \cite{ddj}, \cite{11}, \cite{12}, \cite{13}, \cite{kerstin1}
and \cite{14}) considered the zero location of special classes of $_{2}%
F_{1}\left(  -n,b;c;z\right)  $ with restrictions on the parameters $b$ and
$c$. Results on the asymptotic zero distribution of certain classes of $_{2}%
F_{1}\left(  -n,b;c;z\right)  $ have also appeared (cf. \cite{10}, \cite{sj2},
\cite{15}, \cite{16} and \cite{zhou}).

Different types of $_{2}F_{1}\left(  -n,b;c;z\right)  $ have
well-established connections with classical orthogonal polynomials, notably
the Jacobi polynomials and the Gegenbauer or ultraspherical polynomials (cf. \cite{aar}). For
the ranges of the parameters where these polynomials are orthogonal,
information about the zeros of ${}_{2}F_{1}\left(  -n,b;c;z\right)  $ follows
immediately from classical results (cf. \cite{aar}, \cite{szego}). The
asymptotic zero distribution of $_{2}F_{1}\left(  -n,b;c;z\right)  $ when $b$
and $c$ depend on $n$ can be deduced from recent results by Kuijlaars,
Mart\'{\i}nez-Finkelshtein, Mart\'{\i}nez-Gonz\'{a}lez and Orive (cf.
\cite{kuilmf}, \cite{kuilmfo}, \cite{mfmgo}, \cite{mfo}) on the asymptotic
zero distribution of Jacobi polynomials $P_{n}^{(\alpha,\beta)}(x)$ when the
parameters $\alpha$ and $\beta$ depend on $n$. Conversely, if the distribution
of the zeros of $_{2}F_{1}\left(  -n,b;c;z\right)  $ is known, this leads to
information about the zero distribution of other special functions (cf.
\cite{11}). This makes knowledge of the zero distribution of $_{2}F_{1}\left(
-n,b;c;z\right)  $ extremely valuable.

\medskip The orthogonality of the polynomials $_{2}F_{1}\left(
-n,b;c;z\right)  $ given in the next theorem follows from the orthogonality of
the Jacobi polynomials (cf. \cite[p. 257-261]{rain}) and can also be proved
directly using the Rodrigues' formula for the polynomials ${}_{2}F_{1}\left(
-n,b;c;z\right)  $ (cf. \cite[p. 99]{aar}) as was done in \cite{sj1} and \cite{kuilmfo}.

\begin{theorem}
[cf. \cite{sj1}]\label{basic1} Let $n\in\mathbb{N}_{0}$, $b,\,c\in\mathbb{R}$
and $-c\notin\mathbb{N}_{0}$. Then ${}_{2}F_{1}\left(  -n,b;c;z\right)$ is the
$n^{\textrm{th}}$ degree orthogonal polynomial for the $n$-dependent positive
weight function $|z^{c-1}(1-z)^{b-c-n}|$ on the intervals

\begin{enumerate}
\item[(i)] $(-\infty,0)$ for $c>0$ and $b<1-n$;

\item[(ii)] $(0,1)$ for $c>0$ and $b>c+n-1$;

\item[(iii)] $(1,\infty)$ for $c+n-1<b<1-n$.
\end{enumerate}
\end{theorem}

\begin{figure}[tbh]
\begin{center}
\begin{pspicture}(-4,-3.75)(3.5,3.75)
\psaxes[linewidth=1pt,labels=none,ticks=none]{->}(0,0)(-5,-3.5)(3.5,3.75)
\rput[rt](1,3){\boldmath{$\mathcal{G}_1$}}\rput[rt](2,-2.5){\boldmath{$\mathcal{G}_2$}}
\rput[rt](-4.3,-2.3){\boldmath{$\mathcal{G}_3$}}
\pspolygon[linecolor=white,fillstyle=vlines](0.03,1.15)(2.5,3.5)(0.03,3.5)
\pspolygon[linecolor=white,fillstyle=vlines](0.03,-2)(0.03,-3.5)(3.5,-3.5)(3.5,-2)
\pspolygon[linecolor=white,fillstyle=vlines](-3.4,-2)(-5,-2)(-5,-3.5)
\psline(-5,-3.5)(2.5,3.5)
\psline(-5,-2)(3.5,-2)
\rput(3.5,-0.3){$c$}
\rput(-0.3,3.8){$b$}
\rput(2.5,-1.8){\footnotesize $b=1-n$}
\rput(1.8,2){\footnotesize $b=c+n-1$}
\end{pspicture}
\end{center}
\caption{Values of $b$ and $c$ for which ${}_{2}F_{1}(-n,b;c;z)$ is orthogonal
and has $n$ real simple zeros in the intervals $(0,1)$, $(-\infty,0)$ and
$(1,\infty)$ are indicated by regions $\mathcal{G}_{1}$, $\mathcal{G}_{2}$ and
$\mathcal{G}_{3}$ respectively .}%
\label{orthog1}%
\end{figure}
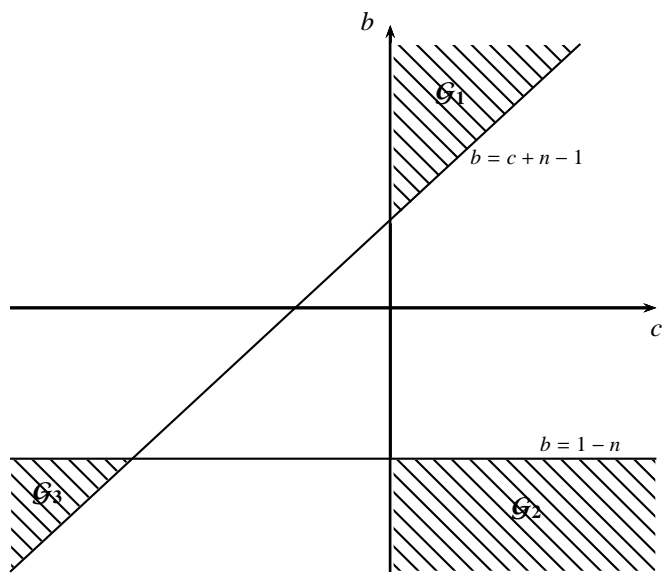

As a consequence of orthogonality, we know that for each $n$, the $n$ zeros of
$_{2}F_{1}\left(  -n,b;c;z\right)  $ are real, simple and lie in the interval
of orthogonality for the corresponding ranges of the parameters (see, for example,
\cite{kerstin2}, Theorem 4) as illustrated in Figure \ref{orthog1}.

\medskip

In his classical paper (cf. \cite{klein}), Felix Klein obtained results on the
precise number of zeros of ${}_{2}F_{1}\left(  a,b;c;z\right)  $ that lie in
each of the intervals $(-\infty,0)$, $(0,1)$ and $(1,\infty)$ by generalizing
earlier results of Hilbert (cf. \cite{Hilbert}). These Hilbert-Klein formulas
are valid for hypergeometric functions and not only for polynomials. Szeg\"{o}
recaptured these results for the special case of Jacobi polynomials
$P_{n}^{(\alpha,\beta)}(x)$, which have a representation as ${}_{2}%
F_{1}\left(  -n,b;c;z\right)  $, in the intervals $(-\infty,-1)$, $(-1,1)$ and
$(1,\infty)$ (cf. \cite{szego}, p.145, Theorem 6.72). The number and location
of the real zeros of ${}_{2}F_{1}\left(  -n,b;c;z\right)  $ for $b$ and $c$
real can be deduced as follows.

\begin{theorem}
[cf. \cite{kerstin1}, Theorem 3.2]\label{basic2} Let $n\in\mathbb{N}$,
$b,\,c\in\mathbb{R}$ and $c>0$. Then,

\begin{itemize}
\item[(i)] For $b>c+n$, all zeros of $_{2}F_{1}\left(  -n,b;c;z\right)  $ are
real and lie in the interval $(0,1)$.

\item[(ii)] For $c<b<c+n$, $c+j-1<b<c+j$, $j=1,2,\dots,n$, $_{2}F_{1}\left(
-n,b;c;z\right)  $ has $j$ real zeros in $(0,1)$. The remaining $(n-j)$ zeros
of $_{2}F_{1}\left(  -n,b;c;z\right)  $ are all non-real if $(n-j)$ is even,
while if $(n-j)$ is odd, $_{2}F_{1}\left(  -n,b;c;z\right)  $ has $(n-j-1)$
non-real zeros and one additional real zero in $(1,\infty)$.

\item[(iii)] For $0<b<c$, all the zeros of $_{2}F_{1}\left(  -n,b;c;z\right)
$ are non-real if $n$ is even, while if $n$ is odd, $_{2}F_{1}\left(
-n,b;c;z\right)  $ has one real zero in $(1,\infty)$ and the other $(n-1)$
zeros are non-real.

\item[(iv)] For $-n<b<0$, $-j<b<-j+1$, $j=1,2,\dots,n$, $_{2}F_{1}\left(
-n,b;c;z\right)  $ has $j$ real negative zeros. The remaining $(n-j)$ zeros of
$_{2}F_{1}\left(  -n,b;c;z\right)  $ are all non-real if $(n-j)$ is even,
while if $(n-j)$ is odd, $_{2}F_{1}\left(  -n,b;c;z\right)  $ has $(n-j-1)$
non-real zeros and one additional real zero in $(1,\infty)$.

\item[(v)] For $b<-n$, all zeros of $_{2}F_{1}\left(  -n,b;c;z\right)  $ are
real and negative.
\end{itemize}
\end{theorem}

The values of the parameters $b$ and $c$ for which ${}_{2}F_{1}\left(
-n,b;c;z\right)  $ has exactly $n$ real simple zeros in $(0,1)$ given in
Theorem \ref{basic2} (i) and (ii) correspond to those in Theorem
\ref{basic1}(ii) while the parameter values in Theorem \ref{basic1}(i) that
ensure that all the zeros of $_{2}F_{1}\left(  -n,b;c;z\right)  $ are real,
simple and negative are the same as those in Theorem \ref{basic2} (iv) and
(v). The values of $b$ and $c$ in Theorem \ref{basic1} (iii) for which $n$
zeros are in $(1,\infty)$ can also be obtained from Theorem \ref{basic2} (iv)
and (v) using the transformation (cf. \cite[p. 79, (2.3.14)]{aar})
\[
{}_{2}F_{1}\left(  -n,b;c;z\right)  =\frac{(c-b)_{n}}{(c)_{n}}{}_{2}%
F_{1}\left(  -n,b;1-n+b-c;1-z\right)
\]
due to Pfaff.

A natural question to ask is whether the parameter ranges in Theorems
\ref{basic1} and \ref{basic2} are the only values of $b$, $c\in\mathbb{R}$ for
which $_{2}F_{1}\left(  -n,b;c;z\right)  $ have $n$ real simple zeros. In this
paper, we use a method that does not rely on orthogonality to determine all
the real values of the parameters $b$ and $c$ for which $_{2}F_{1}\left(
-n,b;c;z\right)  $ have $n$ real simple zeros. We apply an algorithm which
counts the zeros of polynomials with real coefficients and their multiplicities.
We also determine the intervals where the real zeros are located for these
values of $b$ and $c$.

\section{The algorithm}

Recall that given two polynomials $f(x)$ and $g(x)$, with $\deg(f)\geq\deg(g)
$, there exist unique polynomials $q(x)$ and $r (x)$ such that $f (x) =
q(x)g(x) + r (x)$ with $\deg(r) < \deg(g)$. We will denote the leading
coefficient of a polynomial $f(x)=a_{n}x^{n}+a_{n-1}x^{n-1} +\dots+ a_{0}$ by
$lc(f) = a_{n}$.

We use the following algorithm (cf. \cite{rah}).

Let $f(x)$ be a real polynomial with deg$(f)=n\geq2$. Define
\[
f_{0}(x):=f(x)\qquad\text{and}\qquad f_{1}(x):=f^{\prime}(x)
\]
and proceed for $k\in\mathbb{N}$ as follows: 

If deg$(f_{k})>0$ perform the division of $f_{k-1}$ by $f_{k}$ to obtain
\[
f_{k-1}(x)=q_{k-1}(x)f_{k}(x)-r_{k}(x).
\]
Define
\[
f_{k+1}(x)=%
\begin{cases}
r_{k}(x) & \text{if }r_{k}(x)\not \equiv 0\\
f_{k}^{\prime}(x) & \text{if }r_{k}(x)\equiv0
\end{cases}
\]
and generate the sequence of numbers $c_{1},c_{2},\ldots$ where
\[
c_{k}=%
\begin{cases}
\displaystyle\frac{lc(f_{k+1})}{lc(f_{k-1})} & \text{if }r_{k}(x)\not \equiv
0\\
0 & \text{if }r_{k}(x)\equiv0
\end{cases}
.
\]
When $f_{k}$ is constant, the algorithm terminates.

Note that the algorithm must terminate, since the degrees of the polynomials
$f_{k}(x)$ decrease on each step.

\medskip Then we have the following theorem which we will apply to $_{2}%
F_{1}\left(  -n,b;c;z\right)  .$

\begin{theorem}
[cf. \cite{rah}, Theorem 10.5.7, p.339]\label{algorithm2} Let $f$ be a
polynomial of degree $n$ with real coefficients. Then $f$ has only real zeros
if and only if the above algorithm produced $n-1$ non-negative numbers
$c_{1},\ldots,c_{n-1}$. Moreover, the zeros of $f$ are all real and simple if
and only if the numbers $c_{1}, \ldots,c_{n-1}$ are all positive.
\end{theorem}

\section{Main results}

We shall assume throughout our discussion that $b,c\in\mathbb{R}$ with
$b,c\neq0,-1,\dots,-n+1$. The assumption on $b$ is made to ensure that
$_{2}F_{1}\left(  -n,b;c;z\right)  $ is a polynomial of degree $n$.

\begin{proposition}
\label{prop1} Let $b,c\in\mathbb{R}$. Then,

\begin{enumerate}
\item The zeros of $_{2}F_{1}\left(  -2,b;c;z\right)  $ are real and simple if
and only if either (see Figure \ref{fig1}):

\begin{enumerate}
\item[(i)] $c<-1$ and $c<b<0.$

\item[(ii)] $-1<c<0$ and $b>0$ or $b<c.$

\item[(iii)] $c>0$ \ and $b<0$ or $c<b.$
\end{enumerate}

\item The zeros of $_{2}F_{1}\left(  -3,b;c;z\right)  $ are real and simple if
and only if either (see Figure \ref{fig2}):

\begin{enumerate}
\item[(i)] $c<-2$ \ and $1+c<b<-1.$

\item[(ii)] $-2<c<-1$ \ and $-1<b<1+c.$

\item[(iii)] $c>-1,c\neq0$ and $b<-1$ or $b>c+1.$
\end{enumerate}
\end{enumerate}
\end{proposition}

\medskip\begin{figure}[tbh]
\psset{unit=0.75cm}
\par
\begin{center}
\parbox{6cm}{\begin{pspicture}(-3.2,-3.2)(3.2,3.2)
\pspolygon[linecolor=white,fillstyle=vlines](0,0)(0,3)(3,3)
\pspolygon[linecolor=white,fillstyle=vlines](-0,-0)(-0.5,-0.5)(-0.5,-3)(-0,-3)
\pspolygon[linecolor=white,fillstyle=vlines](-0.5,0)(-0.5,3)(-0,3)(-0,0)
\pspolygon[linecolor=white,fillstyle=vlines](-3,-3)(-3,-0)(-0.5,-0)(-0.5,-0.5)
\pspolygon[linecolor=white,fillstyle=vlines](0,-0)(3,-0)(3,-3)(0,-3)
\psaxes[linewidth=1pt,labels=none,ticks=none]{<->}(0,0)(3.2,3.2)(-3.2,-3.2)
\psline[linestyle=dashed](-3,-3)(3,3)
\psline[linestyle=dashed](-0.5,-3)(-0.5,3)
\rput(3.3,-0.1){$c$}
\rput(-0.1,3.3){$b$}
\rput(-0.6,0.1){\tiny{-1}}
\rput(1.2,0.7){\footnotesize{$b=c$}}
\end{pspicture}
\caption{Values of $b$ and $c$ for which ${}_2F_1(-2,b;c;z)$ has only real simple zeros}
\label{fig1}} \qquad\parbox{6cm}{
\begin{pspicture}(-3.1,-3.1)(3.1,3.1)
\pspolygon[linecolor=white,fillstyle=vlines](-0.5,0)(-0.5,3)(2.5,3)
\pspolygon[linecolor=white,fillstyle=vlines](-0.5,0)(-0.5,-0.5)(-1,-0.5)
\pspolygon[linecolor=white,fillstyle=vlines](-0.5,-0.5)(3,-0.5)(3,-3)(-0.5,-3)
\pspolygon[linecolor=white,fillstyle=vlines](-1,-0.5)(-3,-2.5)(-3,-0.5)
\psaxes[linewidth=1pt,labels=none,ticks=none]{<->}(0,0)(-3.2,-3.2)(3.2,3.2)
\psline[linestyle=dashed](-0.5,-3)(-0.5,3)
\psline[linestyle=dashed](-1,-3)(-1,3)
\psline[linestyle=dashed](-3,-0.5)(3,-0.5)
\psline[linestyle=dashed](-3,-2.5)(2.5,3)
\rput(3.3,-0.1){$c$}
\rput(-0.1,3.3){$b$}
\rput(-0.6,0.1){\tiny{-1}}
\rput(-1.1,0.1){\tiny{-2}}
\rput(0.1,-0.4){\tiny{-1}}
\rput(1.3,1.1){\footnotesize{$b=c+1$}}
\end{pspicture}
\caption{Values of $b$ and $c$ for which ${}_2F_1(-3,b;c;z)$ has only real simple zeros}
\label{fig2}}
\end{center}
\end{figure}

\begin{theorem}
\label{main} For any integer $n\geq4$, the polynomial ${}_{2}F_{1}\left(  -n,b;c;z\right)
$ has only real and simple zeros if and only if $\left(  c,b\right)  $ belongs
to one of the four $n$-dependent regions $\mathcal{R}_{1},\ldots
,\mathcal{R}_{4}$ defined by
\begin{align*}
\mathcal{R}_{1} &  =\left\{  c+n-2<b<2-n\right\}  ,\\
\mathcal{R}_{2} &  =\left\{  c>-1,\quad b<2-n\right\}  ,\\
\mathcal{R}_{3} &  =\left\{  c>-1,\quad b>n-2,\quad b>c+n-2\right\}  ,\\
\mathcal{R}_{4} &  =\left\{  -1<c<0,\quad c+n-2<b<n-2\right\}  .
\end{align*}

\end{theorem}

\begin{figure}[tbh]
\begin{center}
\psset{unit=0.5cm} \begin{pspicture}(-8,-7.5)(7,7.5)
\psline{->}(-8,0)(7,0)
\rput(0.2,-0.2){\tiny $0$}
\psline[linewidth=1.2pt,linestyle=dashed]{->}(0,-7.4)(0,7)
\pspolygon[linecolor=white,fillstyle=vlines](0,2)(5,7)(0,7)
\pspolygon[fillstyle=solid,fillcolor=gray](0,2)(0,1.5)(5.5,7)(5,7)
\pspolygon[fillstyle=solid,fillcolor=gray](0,2)(-0.4,1.6)(-0.4,7)(0,7)
\pspolygon[fillstyle=solid,fillcolor=gray](0,2)(0,1.5)(-0.4,1.1)(-0.4,1.6)
\pspolygon[linecolor=white,fillstyle=vlines](0,-4)(0,-7.5)(7,-7.5)(7,-4)
\pspolygon[fillstyle=solid,fillcolor=gray](0,-4)(0,-3.6)(7,-3.6)(7,-4)
\pspolygon[fillstyle=solid,fillcolor=gray](0,-4)(-0.4,-4)(-0.4,-7.5)(0,-7.5)
\pspolygon[fillstyle=solid,fillcolor=gray](0,-4)(-0.4,-4)(-0.4,-3.6)(0,-3.6)
\pspolygon[linecolor=white,fillstyle=vlines](-6,-4)(-8,-4)(-8,-5.7)
\pspolygon[fillstyle=solid,fillcolor=gray](-6,-4)(-5.6,-3.6)(-8,-3.6)(-8,-4)
\pspolygon[fillstyle=solid,fillcolor=gray](-6,-4)(-5.5,-4)(-8,-6.5)(-8,-6)
\pspolygon[fillstyle=solid,fillcolor=gray](-6,-4)(-5.5,-4)(-5.1,-3.6)(-5.6,-3.6)
\psline[linestyle=dashed](-8,-6.5)(5.5,7)
\psline[linestyle=dashed](-0.4,-7.4)(-0.4,7)
\psline[linestyle=dashed](-8,-6)(5,7)
\psline[linestyle=dashed](-8,-4)(7,-4)
\psline[linestyle=dashed](-8,-3.6)(7,-3.6)
\rput(7,-0.3){$c$}
\rput(-0.3,7.3){$b$}
\rput(-0.65,-0.2){\tiny $-1$}
\rput(-1.3,6){\footnotesize $c=-1$}
\rput(-2,-2.3){\footnotesize $b=c+n-2$}
\rput(-5.3,-1.5){\footnotesize $b=c+n-1$}
\rput(5,-3.2){\footnotesize $b=2-n$}
\rput(8,-4.4){\footnotesize $b=1-n$}
\rput(0.9,-1.5){\footnotesize $c=0$}
\psline[linestyle=dashed](-8,1.5)(7,1.5)
\rput(5,1.9){\footnotesize $b=n-2$}
\end{pspicture}
\end{center}
\caption{Values of $b$ and $c$ for which ${}_{2}F_{1}(-n,b;c;z)$,
$n=4,5,\dots$ has $n$ real simple zeros}%
\label{all}%
\end{figure}
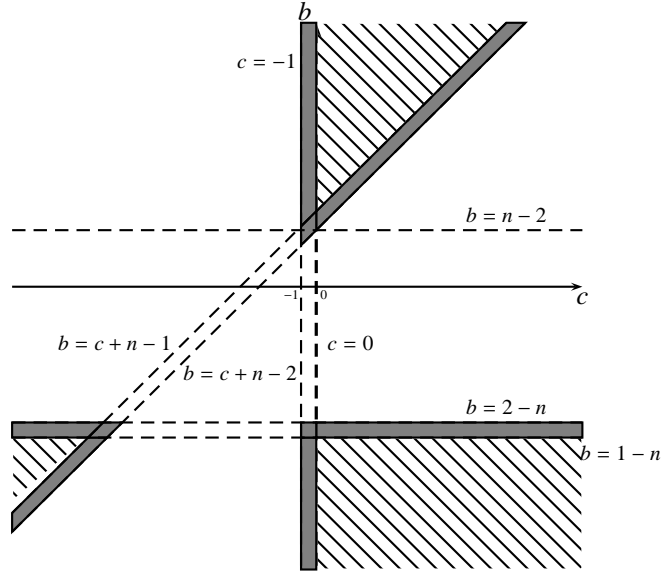

The parameter values $\left(  c,b\right)  \in\mathcal{R}_{1}\cup
\mathcal{R}_{2}\cup\mathcal{R}_{3}\cup\mathcal{R}_{4}$ described in Theorem
\ref{main} for which ${}_{2}F_{1}\left(  -n,b;c;z\right)  $, $n=4,5,\dots$ has
$n$ real simple zeros are illustrated by the grey and diagonally shaded
regions in Figure \ref{all} with the grey regions indicating those parameter
values that extend the results in Theorems \ref{basic1} and \ref{basic2}.

Next, we turn our attention to the location of the zeros of $_{2}F_{1}\left(
-n,b;c;z\right)  $ for those parameter values located in the grey shaded
regions in Figure \ref{all} where the polynomials are no longer orthogonal and
the location of the real zeros cannot be obtained using Theorems \ref{basic1}
and \ref{basic2}. For these values of $b$ and $c$, the polynomials $_{2}%
F_{1}\left(  -n,b;c;z\right)  $ are quasi-orthogonal of order 1 and, in some
cases, order 2 (cf. \cite{chihara} and \cite{Brez}). Theorem 3 in \cite{Brez}
and Theorem 6 in \cite{Joulak} yield information on the zero location of
quasi-orthogonal polynomials with non-varying weight functions. However, these
results cannot be applied to $_{2}F_{1}\left(  -n,b;c;z\right)  $ since their
weight function depends on $n$. We use information about the zeros of Jacobi
polynomials to obtain the following three results.

\begin{theorem}
\label{jacobi2f1zeros1} Let $n\in\mathbb{N}$ and $b,c\in\mathbb{R}$. Then,
$_{2}F_{1}\left(  -n,b;c;z\right)  $ has all its zeros real and simple and
$(n-2)$ of them lie in

\begin{enumerate}
\item[(i)] $(0,1)$ for $-1<c<0$ and $c+n-2<b<c+n-1$. One of the remaining
zeros lies in $(1,\infty)$ and the other one in $(-\infty,0)$.

\item[(ii)] $(1,\infty)$ for $1-n<b<2-n$ and $b-n+1<c<b-n+2$. One of the
remaining zeros lies in $(-\infty,0)$ and the other one in $(0,1)$.

\item[(iii)] $(-\infty,0)$ for $-1<c<0$ and $1-n<b<2-n$. One of the remaining
zeros lies in $(1,\infty)$ and the other one in $(0,1)$.
\end{enumerate}
\end{theorem}

Theorem \ref{jacobi2f1zeros1} applies to the parameter values illustrated in
Figure \ref{qo2}.

\begin{figure}[tbh]
\begin{center}
\psset{unit=0.5cm} \begin{pspicture}(-8,-7.5)(7,7.5)
\pspolygon[fillstyle=solid,fillcolor=gray](0,2)(0,1.5)(-0.4,1.1)(-0.4,1.6)
\pspolygon[fillstyle=solid,fillcolor=gray](0,-4)(-0.4,-4)(-0.4,-3.6)(0,-3.6)
\pspolygon[fillstyle=solid,fillcolor=gray](-6,-4)(-5.5,-4)(-5.1,-3.6)(-5.6,-3.6)
\psline{->}(-8,0)(7,0)
\rput(0.2,-0.2){\tiny $0$}
\psline[linewidth=1.2pt,linestyle=dashed]{->}(0,-7.4)(0,7)
\psline[linestyle=dashed](-8,-6.5)(5.5,7)
\psline[linestyle=dashed](-0.4,-7.4)(-0.4,7)
\psline[linestyle=dashed](-8,-6)(5,7)
\psline[linestyle=dashed](-8,-4)(7,-4)
\psline[linestyle=dashed](-8,-3.6)(7,-3.6)
\rput(7,-0.3){$c$}
\rput(-0.3,7.3){$b$}
\rput(-0.6,-0.2){\tiny $-1$}
\rput(-1.3,6){\footnotesize $c=-1$}
\rput(-2,-2.3){\footnotesize $b=c+n-2$}
\rput(-5.3,-1.5){\footnotesize $b=c+n-1$}
\rput(5,-3.2){\footnotesize $b=2-n$}
\rput(6.5,-4.4){\footnotesize $b=1-n$}
\rput(0.9,-1.5){\footnotesize $c=0$}
\psline[linestyle=dashed](-8,1.5)(7,1.5)
\rput(5,1.9){\footnotesize $b=n-2$}
\end{pspicture}
\end{center}
\caption{Values of $b$ and $c$ corresponding to those described in Theorem
\ref{jacobi2f1zeros1}}%
\label{qo2}%
\end{figure}

\begin{theorem}
\label{jacobi2f1zeros2} Let $n\in\mathbb{N}$ and $b,c\in\mathbb{R}$. Then,
$_{2}F_{1}\left(  -n,b;c;z\right)  $ has all its zeros real and simple and
$(n-1)$ of them lie in

\begin{enumerate}
\item[(i)] $(0,1)$ for $-1<c<0$ and $b>c+n-1$. The remaining zero is negative.

\item[(ii)] $(1,\infty)$ for $1-n<b<2-n$ and $c<b-n-1$. The remaining zero is negative.

\item[(iii)] $(-\infty,0)$ for $-1<c<0$ and $b<1-n$. The remaining zero lies
in $(0,1)$.
\end{enumerate}
\end{theorem}

The parameter values described in Theorem \ref{jacobi2f1zeros2} are
illustrated in Figure \ref{qo1a}.

\begin{figure}[tbh]
\begin{center}
\psset{unit=0.5cm} \begin{pspicture}(-8,-7.5)(7,7.5)
\pspolygon[fillstyle=solid,fillcolor=gray](0,2)(-0.4,1.6)(-0.4,7)(0,7)
\pspolygon[fillstyle=solid,fillcolor=gray](0,-4)(-0.4,-4)(-0.4,-7.5)(0,-7.5)
\pspolygon[fillstyle=solid,fillcolor=gray](-6,-4)(-5.6,-3.6)(-8,-3.6)(-8,-4)
\psline{->}(-8,0)(7,0)
\rput(0.2,-0.2){\tiny $0$}
\psline[linewidth=1.2pt,linestyle=dashed]{->}(0,-7.4)(0,7)
\psline[linestyle=dashed](-8,-6.5)(5.5,7)
\psline[linestyle=dashed](-0.4,-7.4)(-0.4,7)
\psline[linestyle=dashed](-8,-6)(5,7)
\psline[linestyle=dashed](-8,-4)(7,-4)
\psline[linestyle=dashed](-8,-3.6)(7,-3.6)
\rput(7,-0.3){$c$}
\rput(-0.3,7.3){$b$}
\rput(-0.6,-0.2){\tiny $-1$}
\rput(-1.3,6){\footnotesize $c=-1$}
\rput(-2,-2.3){\footnotesize $b=c+n-2$}
\rput(-5.3,-1.5){\footnotesize $b=c+n-1$}
\rput(5,-3.2){\footnotesize $b=2-n$}
\rput(8,-4.4){\footnotesize $b=1-n$}
\rput(0.9,-1.5){\footnotesize $c=0$}
\end{pspicture}
\end{center}
\caption{Values of $b$ and $c$ corresponding to those described in Theorem
\ref{jacobi2f1zeros2}}%
\label{qo1a}%
\end{figure}

\begin{theorem}
\label{jacobi2f1zeros3} Let $n\in\mathbb{N}$ and $b,c\in\mathbb{R}$. Then,
$_{2}F_{1}\left(  -n,b;c;z\right)  $ has all its zeros real and simple and
$(n-1)$ of them lie in

\begin{enumerate}
\item[(i)] $(0,1)$ for $c>0$ and $c+n-2<b<c+n-1$. The remaining zero is in the
interval $(1,\infty)$.

\item[(ii)] $(1,\infty)$ for $b<1-n$ and $c+n-2<b<c+n-1$. The remaining zero
lies in $(0,1)$.

\item[(iii)] $(-\infty,0)$ for $c>0$ and $1-n<b<2-n$. The remaining zero lies
in $(1,\infty)$.
\end{enumerate}
\end{theorem}

Figure \ref{qo1b} illustrates the range of the parameters $b$ and $c$ referred
to in Theorem \ref{jacobi2f1zeros3}.

\begin{figure}[tbh]
\begin{center}
\psset{unit=0.5cm} \begin{pspicture}(-8,-7.5)(7,7.5)
\pspolygon[fillstyle=solid,fillcolor=gray](0,2)(0,1.5)(5.5,7)(5,7)
\pspolygon[fillstyle=solid,fillcolor=gray](0,-4)(0,-3.6)(7,-3.6)(7,-4)
\pspolygon[fillstyle=solid,fillcolor=gray](-6,-4)(-5.5,-4)(-8,-6.5)(-8,-6)
\psline{->}(-8,0)(7,0)
\rput(0.2,-0.2){\tiny $0$}
\psline[linewidth=1.2pt,linestyle=dashed]{->}(0,-7.4)(0,7)
\psline[linestyle=dashed](-8,-6.5)(5.5,7)
\psline[linestyle=dashed](-0.4,-7.4)(-0.4,7)
\psline[linestyle=dashed](-8,-6)(5,7)
\psline[linestyle=dashed](-8,-4)(7,-4)
\psline[linestyle=dashed](-8,-3.6)(7,-3.6)
\rput(7,-0.3){$c$}
\rput(-0.3,7.3){$b$}
\rput(-0.6,-0.2){\tiny $-1$}
\rput(-1.3,6){\footnotesize $c=-1$}
\rput(-2,-2.3){\footnotesize $b=c+n-2$}
\rput(-5.3,-1.5){\footnotesize $b=c+n-1$}
\rput(5,-3.2){\footnotesize $b=2-n$}
\rput(8,-4.4){\footnotesize $b=1-n$}
\rput(0.9,-1.5){\footnotesize $c=0$}
\end{pspicture}
\end{center}
\caption{Values of $b$ and $c$ corresponding to those described in Theorem
\ref{jacobi2f1zeros3}}%
\label{qo1b}%
\end{figure}

\section{Proofs}

\proof{Proposition \ref{prop1}}.

\begin{enumerate}
\item Since
\[
_{2}F_{1}\left(  -2,b;c;z\right)  =1-\frac{2b}{c}z+\frac{b\left(  b+1\right)
}{c\left(  c+1\right)  }z^{2},
\]
we see that $_{2}F_{1}\left(  -2,b;c;z\right)  =0$ if and only if
\[
z=\frac{b\left(  c+1\right)  \pm\sqrt{b\left(  c+1\right)  \left(  b-c\right)
}}{b\left(  b+1\right)  }.
\]
Hence, the zeros of $_{2}F_{1}\left(  -2,b;c;z\right)  $ are real and simple
if and only if $b\left(  c+1\right)  \left(  b-c\right)  >0.$

\item The discriminant of
\[
_{2}F_{1}\left(  -3,b;c;z\right)  =1-\frac{3b}{c}z+\frac{3b\left(  b+1\right)
}{c\left(  c+1\right)  }z^{2}-\frac{b\left(  b+1\right)  \left(  b+2\right)
}{c\left(  c+1\right)  \left(  c+2\right)  }z^{3}%
\]
is given by $\displaystyle{\Delta_{3}=108\frac{b^{2}(b+1)(b-c-1)(b-c)^{2}%
}{c^{4}(c+1)^{3}(c+2)^{2}}}$ (cf. \cite{Gelfond}) and therefore $_{2}%
F_{1}\left(  -3,b;c;z\right)  $ has real simple roots if and only if
$\Delta_{3}>0$.
\end{enumerate}

\hfill\hspace{-6pt}\rule[-4pt]{6pt}{6pt} \vskip8pt plus3pt minus 3pt

The following two lemmas will be used in the proof of our main result.

\begin{lemma}
\label{Lemma 1} Let
\[
\alpha_{k,l}=\frac{\left(  \frac{k-n}{2}\right)  _{l}\left(  \frac
{2k-n-b-3}{4}\right)  _{l}\left(  \frac{k-b-1+c}{2}\right)  _{l}}{\left(
\frac{k-b-2}{2}\right)  _{l}\left(  \frac{2k-b-1-n}{4}\right)  _{l}\left(
\frac{k-n-1-c}{2}\right)  _{l}}%
\]
and let the sequence $\theta_{k}$ be recursively defined by
\begin{equation}
\theta_{k+1}=\alpha_{k,1}~\theta_{k-1},\label{reqTheta}%
\end{equation}
for $k\in\{2,\dots,n-2\}$, with
\[
\theta_{1}=-\frac{nb}{c},\quad\theta_{2}=\frac{\left(  b-c\right)  \left(
n-1\right)  }{c\left(  b+n-1\right)  }.
\]
Then,
\begin{equation}
\theta_{2k}=\frac{(b+1)(n+c)}{c(n+b+1)}\alpha_{1,k}~,\quad k=1,2,\ldots
\left\lfloor \frac{n-1}{2}\right\rfloor ,\label{2ktheta}%
\end{equation}
and
\begin{equation}
\theta_{2k+1}=-\frac{nb}{c}\alpha_{2,k}~,\quad k=0,1,\ldots\left\lfloor
\frac{n-2}{2}\right\rfloor .\label{2kplus1theta}%
\end{equation}

\end{lemma}

\proof{Lemma \ref{Lemma 1}}. We prove the result by induction on $k$. When
$k=1$, the right-hand side of \eqref{2ktheta} is
\[
\frac{(b+1)(n+c) \left(  \frac{1-n}{2} \right)  \left(  -\frac{n+b+1}%
{4}\right)  \left(  \frac{c-b}{2}\right)  }{c(n+b+1)\left(  - \frac{b+1}%
{2}\right)  \left(  -\frac{n+b-1}{4}\right)  \left(  - \frac{n+c}{2}\right)  }
= \frac{(b-c)(n-1)}{c(b+n-1)}%
\]
which is $\theta_{2}$ as required.\newline

We now assume the result is true for $k=t$ and prove the result true for
$k=t+1$. If we let $k=t+1$ on the right-hand side of \eqref{2ktheta}, we
obtain
\begin{align*}
\text{RHS }  &  = \frac{(b+1)(n+c)}{c(n+b+1)}~\alpha_{1,t+1}\\
&  =\frac{(b+1)(n+c)}{c(n+b+1)} \frac{\left(  \frac{1-n}{2}+t\right)  \left(
\frac{c-b}{2}+t\right)  \left(  -\frac{n+b+1}{4}+t \right)  }{\left(
-\frac{b+1}{2}+t\right)  \left(  -\frac{n+c}{2}+t\right)  \left(
-\frac{n+b-1}{4}+t\right)  }~\alpha_{1,t}~~\mbox{since}~~(a)_{k+1}%
=(a+k)(a)_{k}\\
&  = \theta_{2t} \frac{\left(  1-n+2t\right)  \left(  c-b+2t\right)  \left(
-n-b-1+4t \right)  }{\left(  -b-1+2t \right)  \left(  -n-c+2t\right)  \left(
-n-b+1+4t\right)  }~~\mbox{by the inductive hypothesis}\\
&  = \alpha_{2t+1,1}~\theta_{2t}\\
&  = \theta_{2t+2} ~~\mbox{from (\ref{reqTheta})}
\end{align*}
and the result follows by induction.

\medskip The second relation \eqref{2kplus1theta} may be proved by induction
in a similar way. \hfill\hspace{-6pt}\rule[-4pt]{6pt}{6pt} \vskip8pt plus3pt
minus 3pt

\begin{lemma}
\label{Lemma 2} Let $n\geq4$. Then, for all $k\in\{2,\dots,n-1\}$,
\begin{equation}
\frac{(n-k)(n+c-k)(b+1-k)(b-c+1-k)}{(n+b+2-2k)(n+b-2k)(n+b+1-2k)^{2}%
}\label{ck}%
\end{equation}
and
\begin{equation}
\frac{\left(  n-1\right)  \left(  n+c-1\right)  \left(  b-c\right)  }{\left(
n+b-2\right)  \left(  n+b-1\right)  ^{2}}\label{c1}%
\end{equation}
are positive if and only if $\left(  c,b\right)  \in\mathcal{R}_{1}%
\cup\mathcal{R}_{2}\cup\mathcal{R}_{3}\cup\mathcal{R}_{4}$.
\end{lemma}

\medskip\proof{Lemma \ref{Lemma 2}}. Since $n-1>0$ and $(n+b-1)^{2}>0$ for all
$n\in\mathbb{N}$, $b\in\mathbb{R}$, we see that (\ref{c1}) is positive if and
only if
\begin{align*}
(c,b)  &  \in\{c<1-n,~b>c,~b<2-n\}=\mathcal{A}_{1}\supset\mathcal{R}_{1}~
\text{or}\\
&  \in\{c>1-n,~b<c,~b<2-n\}=\mathcal{A}_{2}\supset\mathcal{R}_{2}~ \text{or}\\
&  \in\{c>1-n,~b>c,~b>2-n\}=\mathcal{A}_{3}\supset(\mathcal{R}_{3}%
\cup\mathcal{R}_{4})~ \text{or}\\
&  \in\{c<1-n,~b<c,~b>2-n\}= \emptyset.
\end{align*}

Clearly $\displaystyle{\frac{(n-k)}{(n+b+1-2k)^{2}}>0}$ for all $k\in
\{2,\dots,n-1\}$, $n\in\mathbb{N}$ and $b\in\mathbb{R}$, $b\neq,
3-n,5-n,\dots,n-5,n-3$.

\medskip Furthermore $b>n-2$ if and only if $n+b+2-2k>0$, $n+b-2k>0$ and
$b+1-k>0$ for all $k\in\{2,\dots,n-1\}$. Hence, when $b>n-2$, (\ref{ck}) will
be positive for all $k\in\{2,\dots,n-1\}$ if and only if
\begin{align*}
(c,b)  &  \in\{b>c+k-1,~c+n-k>0,~k=2,\dots,n-1\}=
\{b>c+n-2,~c>-1\}=\mathcal{R}_{3} ~~ \text{or}\\
&  \in\{b<c+k-1,~c+n-k<0,~k=2,\dots,n-1\}= \{b<c+1,~c<2-n\}=\emptyset~~
\text{or}\\
&  \in\{b<c+k-1,~c+n-k<0,~k=2,\dots,l\}\cap\{b>c+k-1,~c+n-k>0,~k=l+1,\dots
,n-1\}=\emptyset~~ \text{or}\\
&  \in\{b>c+k-1,~c+n-k>0,~k=2,\dots,l\}\cap\\
&  \{b<c+k-1,~c+n-k<0,~k=l+1,\dots,n-1\}\cap\{b>n-2\}=\emptyset.
\end{align*}

\medskip Similarly, $b<2-n$ if and only if $b+1-k<0$, $n+b+2-2k<0$ and
$n+b-2k<0$ for $k\in\{2,\dots,n-1\}$. Hence, when $b<2-n$, (\ref{ck}) will be
positive for all $k\in\{2,\dots,n-1\}$ if and only if
\begin{align*}
(c,b)  &  \in\{b>c+k-1,~c+n-k<0,~k=2,\dots,n-1\}=
\{b>c+n-2,~c<2-n\}=\mathcal{A}_{4}\supset\mathcal{R}_{1} ~~ \text{or}\\
&  \in\{b<c+k-1,~c+n-k>0,~k=2,\dots,n-1\}= \{b<c+1,~c>-1\}=\mathcal{A}%
_{5}\supset\mathcal{R}_{2}~~ \text{or}\\
&  \in\{b>c+k-1,~c+n-k<0,~k=2,\dots,l\}\cap\{b<c+k-1,~c+n-k>0,~k=l+1,\dots
,n-1\}=\emptyset~~ \text{or}\\
&  \in\{b<c+k-1,~c+n-k>0,~k=2,\dots,l\}\cap\{b>c+k-1,~c+n-k<0,~k=l+1,\dots
,n-1\}=\emptyset.
\end{align*}

\medskip For the remaining case where $2-n<b<n-2$ or, more specifically,
$-2<b+n-2k<0$ with $b-k+1<0$ for all $k\in\{2,\dots,n-1\}$, the only non-empty
possibility is that (\ref{ck}) is positive for $k=n-1$ if and only if $c>-1$,
$b>c+n-2$ and $n-4<b<n-2$ whereas (\ref{ck}) is positive for $k\in
\{2,3,\dots,n-2\}$ if and only if $c>-1$, $b>c+n-2$ and $b>n-3$. Hence, when
$2-n<b<n-2$, (\ref{ck}) is positive for all $k\in\{2,\dots,n-1\}$ if and only
if $(c,b)\in\mathcal{R}_{4}$.

\medskip Since $\mathcal{A}_{1}\cap\mathcal{A}_{4} = \mathcal{R}_{1}$ and
$\mathcal{A}_{2} \cap\mathcal{A}_{5} = \mathcal{R}_{2}$, the result follows.
\hfill\hspace{-6pt}\rule[-4pt]{6pt}{6pt} \vskip8pt plus3pt minus 3pt

\proof{Theorem \ref{main}}. We apply the algorithm to the polynomial
\[
f(z)=\ _{2}F_{1}\left(  -n,b;c;z\right)  .
\]
We have (cf. \cite{rain}, p.69, ex.1)
\[
f_{1}(z)=f^{\prime}(z)=-\frac{nb}{c}\ _{2}F_{1}\left(  -n+1,b+1;c+1;z\right)
.
\]
Using Raimundas Vid\~{u}nas' Maple package for contiguous relations of
$_{2}F_{1}$ hypergeometric series (cf. \cite{vin}, \cite{vin2}), we obtain
\[
f_{0}(z)=\frac{1}{n}\left(  z-\frac{c+n-1}{b+n-1}\right)  \ f_{1}%
(z)-\frac{\left(  b-c\right)  \left(  n-1\right)  }{c\left(  b+n-1\right)
}\ _{2}F_{1}\left(  -n+2,b;c+1;z\right)  .
\]
This relation can easily be verified by comparing coefficients.
Thus,
\[
\ f_{2}(z)=r_{1}(z)=\frac{\left(  b-c\right)  \left(  n-1\right)  }{c\left(
b+n-1\right)  }\ _{2}F_{1}\left(  -n+2,b;c+1;z\right)  .
\]

In the next step ($k=2)$, we get
\[
f_{1}(z)=q_{1}(z)f_{2}(z)-r_{2}(z),
\]
with
\[
q_{1}\left(  z\right)  =\frac{n\left(  b+n-1\right)  ^{2}\left(  b+n-2\right)
} {\left(  n-1\right)  \left(  c+n-1\right)  \left(  b-c\right)  } \left[  z
+\frac{\left(  n-2\right)  \left(  c+n-2\right)  } {b+n-3}-\frac{\left(
n-1\right)  \left(  c+n-1\right)  }{b+n-1}\right]
\]
and
\[
r_{2}\left(  z\right)  =\frac{\left(  b+n-1\right)  \left(  b-1-c\right)
n\left(  n-2\right)  }{c\left(  c+n-1\right)  \left(  b+n-3\right)  }
\ _{2}F_{1} \left(  -n+3,b-1;c+1;z\right)  .
\]

Setting
\[
f_{k}(z)=\theta_{k}\ _{2}F_{1}\left(  -n+k,b+2-k;c+1;z\right)  ,\quad
k\in\{1,\dots,n-1\},
\]
we see that in general we need a contiguous relation of the form
\begin{align*}
&  \ _{2}F_{1}\left(  -n+k-1,b+3-k;c+1;z\right)  \\
&  \qquad=q_{k-1}\left(  z\right)  \frac{\theta_{k}}{\theta_{k-1}}\ _{2}%
F_{1}\left(  -n+k,b+2-k;c+1;z\right)  -\frac{\theta_{k+1}}{\theta_{k-1}}%
\ _{2}F_{1}\left(  -n+k+1,b+1-k;c+1;z\right)  ,
\end{align*}
for $k\in\{2,\dots,n-2\},$ with
\[
\theta_{1}=-\frac{nb}{c},\quad\theta_{2}=\frac{\left(  b-c\right)  \left(
n-1\right)  }{c\left(  b+n-1\right)  }.
\]

Using Vid\~{u}nas' package, we obtain (\ref{reqTheta}) for $k=2,3,\dots$ 
and from Lemma \ref{Lemma 1} we conclude that $\theta_{k}$ is well defined and
non-zero for all $k\in\{1,\dots,n-1\}$ when $c\neq0$ and $\left(  c,b\right)
\in\mathcal{R}_{1}\cup\mathcal{R}_{2}\cup\mathcal{R}_{3}\cup\mathcal{R}_{4}.$
Thus,
\[
c_{k}=\dfrac{lc\left(  f_{k+1}\right)  }{lc\left(  f_{k-1}\right)  },\quad
k\in\mathbb{N},
\]
which implies that
\[
c_{1}=\frac{\left(  b-c\right)  \left(  n-1\right)  }{c\left(  b+n-1\right)
}\frac{\left(  -n+2\right)  _{n-2}\left(  b\right)  _{n-2}}{\left(
c+1\right)  _{n-2}\left(  n-2\right)  !}\frac{\left(  c\right)  _{n}%
\ n!}{\left(  -n\right)  _{n}\left(  b\right)  _{n}}=\frac{\left(  n-1\right)
\left(  c+n-1\right)  \left(  b-c\right)  }{\left(  b+n-2\right)  \left(
b+n-1\right)  ^{2}}%
\]
and, for $k\in\{2,\dots,n-1\}$,%

\begin{align*}
c_{k}  &  =\frac{\theta_{k+1}}{\theta_{k-1}}\frac{\left(  -n+k+1\right)
_{n-k-1}\left(  b+1-k\right)  _{n-k-1}}{\left(  c+1\right)  _{n-k-1}\left(
n-k-1\right)  !} \frac{\left(  n-k+1\right)  !\left(  c+1\right)  _{n-k+1}%
}{\left(  -n+k-1\right)  _{n-k+1}\left(  b+3-k\right)  _{n-k+1}}\\
&  =\frac{(n-k)(n-k+c)(b-k+1)(-c+1+b-k)}{(n-2k+b+2)(n-2k+b)(b-2k+1+n)^{2}}.
\end{align*}

From Lemma \ref{Lemma 2}, we know that $c_{k}>0$ for all $k\in\{1,\dots,n-1\}$
when $\left(  c,b\right)  \in\mathcal{R}_{1}\cup\mathcal{R}_{2}\cup
\mathcal{R}_{3}\cup\mathcal{R}_{4}$. The result now follows from Theorem
\ref{algorithm2}.\hfill\hspace{-6pt}\rule[-4pt]{6pt}{6pt} \vskip8pt plus3pt
minus 3pt

\proof{Theorem \ref{jacobi2f1zeros1}}. From \cite{Brez}, Corollary 4 (i), we
know that for $-1<\alpha<0$ and $-1<\beta<0$, the Jacobi polynomials
$P_{n}^{(\alpha-1,\beta-1)}(x)$ have real simple zeros and $(n-2)$ of them are
in the interval $(-1,1)$. The smallest zero is smaller than $-1$ and the
largest zero is larger than $1$. Equivalently, the same is true for the zeros
of Jacobi polynomials $P_{n}^{(\alpha,\beta)}(x)$ when $-2<\alpha<-1 $ and
$-2<\beta<-1$.

\begin{enumerate}
\item[(i)] One of the connections between Jacobi polynomials and the
polynomials $_{2}F_{1}\left(  -n,b;c;z\right)  $ is given by (cf. \cite{rain},
p. 254, eq. 3)
\begin{equation}
P_{n}^{(\alpha,\beta)}(x)=\frac{(-1)^{n}(1+\beta)_{n}}{n!}{}_{2}F_{1}\left(
-n,1+\alpha+\beta+n;1+\beta;\frac{x+1}{2}\right)  \label{jachyp1}%
\end{equation}
where $\alpha=b-n-c$ and $\beta=c-1$. The conditions $-2<\beta<-1$ and
$-2<\alpha<-1$ are equivalent to $-1<c<0$ and $c+n-2<b<c+n-1$. Furthermore,
the intervals $(-1,1)$, $(1,\infty)$ and $(-\infty,-1)$ are transformed to
$(0,1)$, $(1,\infty)$ and $(-\infty,-1)$ respectively under the linear mapping
$x=2z-1$. Thus, when $c\in(-1,0)$ and $b\in(c+n-2,c+n-1),$ $_{2}F_{1}\left(
-n,b;c;z\right)  $ has $n-2$ real, simple zeros in the interval $(0,1)$, one
zero in $(1,\infty)$ and one zero in $(-\infty,0)$ for each $n\in\mathbb{N}$.

\item[(ii)] The representation (\cite{rain}, p. 255, eq. 9)
\begin{equation}
P_{n}^{(\alpha,\beta)}(x)=\frac{(1+\alpha+\beta)_{2n}}{n!(1+\alpha+\beta)_{n}%
}\left(  \frac{x+1}{2}\right)  ^{n}{}_{2}F_{1}\left(  -n,-\beta-n;-\alpha
-\beta-2n;\frac{2}{x+1}\right)  \label{jachyp2}%
\end{equation}
where $\alpha=b-c-n$ and $\beta=-b-n$ yields the stated result, since the
restrictions $-2<\alpha<-1$ and $-2<\beta<-1$ will correspond to
$1+b-n<c<2+b-n$ and $1-n<b<2-n$ while the intervals $-1<x<1$, $x>1$ and $x<-1$
are mapped to $z>1$, $0<z<1$ and $z<0$ respectively under the fractional
transformation $\displaystyle{z=\frac{2}{x+1}}$.

\item[(iii)] For this case we use the representation (cf. \cite{rain}, p. 255,
eq. 8)
\begin{equation}
\label{jachyp3}P_{n}^{(\alpha,\beta)}(x) = \frac{(1+\beta)_{n}}{n!} \left(
\frac{x-1}{2}\right)  ^{n} {}_{2}F_{1} \left(  -n, -\alpha-n; \beta+1;
\frac{x+1}{x-1}\right)
\end{equation}
where $\alpha=-b -n$ and $\beta=c-1$. Under the transformation
$\displaystyle{z=\frac{x+1}{x-1}}$, the interval $-1<x<1$ is mapped to the
negative real line while the intervals $x>1$ and $x<-1$ are mapped to
$\displaystyle{z=1+\frac{2}{x-1}>1}$ and $0<z<1$ respectively. Also, since
$-2<\alpha<-1$ and $-2<\beta<-1$ correspond to $1-n<b<2-n$ and $-1<c<0$, the
result follows.
\end{enumerate}

\hfill\hspace{-6pt}\rule[-4pt]{6pt}{6pt} \vskip8pt plus3pt minus 3pt

\proof{Theorem \ref{jacobi2f1zeros2}}. From \cite{Brez}, Corollary 4 (ii) (a),
we know that for $-1<\alpha$ and $-1<\beta<0$, the zeros of the Jacobi
polynomial $P_{n}^{(\alpha,\beta-1)}(x)$ are all real and simple, and $(n-1)$
of them lie in the interval $(-1,1)$. The smallest zero is less than $-1$. We
may equivalently say the same about the zeros of the Jacobi polynomial
$P_{n}^{(\alpha,\beta)}(x)$ for $\alpha>-1$ and $-2<\beta<-1$. Using
transformations \eqref{jachyp1}, \eqref{jachyp2} and \eqref{jachyp3} the
results can be proven following the same reasoning as in the proof of Theorem
\ref{jacobi2f1zeros1}.\hfill\hspace{-6pt}\rule[-4pt]{6pt}{6pt} \vskip8pt
plus3pt minus 3pt

\proof{Theorem \ref{jacobi2f1zeros3}}. From \cite{Brez}, Corollary 4 (ii) (b),
we know that for $-1<\alpha<0$ and $-1<\beta$, the zeros of the Jacobi
polynomial $P_{n}^{(\alpha-1,\beta)}(x)$ are all real and simple, and $(n-1)$
of them lie in the interval $(-1,1)$. Moreover, the largest zero is greater
than $1$. We may equivalently say the same about the zeros of the Jacobi
polynomial $P_{n}^{(\alpha,\beta)}(x)$ for $-2<\alpha<-1$ and $\beta>-1$. The
results now follow from \eqref{jachyp1}, \eqref{jachyp2} and \eqref{jachyp3}
as before.\hfill\hspace{-6pt}\rule[-4pt]{6pt}{6pt} \vskip8pt plus3pt minus 3pt

\end{document}